\input AHTOH-E.STY

\UDC{
512.543.72 
}

\title{The centre of a finitely generated strongly verbally closed group
is almost always pure}

\author{%
Filipp D. Denissov
\quad
Anton A. Klyachko
}
\address{
\myAddressWC
\qquad\qquad
denissov.filipp@gmail.com
}

\grantsSecond{\RSF 22-11-00075}

\abstract{
The assertion in the title
implies that many interesting groups (e.g.,
all non-abelian braid groups or $\SL_{100}(\Z)$)
are not strongly
verbally closed, i.e., they
embed into some
finitely generated groups as verbally closed subgroups, which are not
retracts.
}

\s 1.
Introduction

A subgroup $H$ of a group $G$ is called \emph{verbally closed}
[MR14]
if
any equation of the form
$$
w(x,y,\dots)=h,
\qbox{where $w$ is an element of a free group
$F(x,y,\dots)$ and $h\in H$
},
$$
admits a solution in $G$ if and only if it admits a solution in $H$.
Similarly,
$H$ is
called \emph{algebraically closed} in $G$
if each finite system of equations with coefficients from~$H$
$$
\{w_1(x,y,\dots)=1, \dots, w_m(x,y,\dots)=1\},
\qbox{where $w_i\in H*F(x,y,\dots)$ (and $*$ means the
free product)}
$$
admits a solution in $G$ if and only if it admits a solution in $H$.

The algebraic closedness is a stronger property than the verbal
closedness,
but, often, these properties turn out to be
equivalent (see
[Rom12],
[RKh13],
[MR14],
[Mazh17],
[RKhK17],
[KM18],
[KMM18],
[Mazh18],
[Bog22],
[Mazh19],
[RT20],
[Tim21],
[KMO23],
[KK24]).
A group $H$ is called \emph{strongly verbally closed}
[Mazh18]
if it is algebraically closed
in any group containing $H$ as a verbally closed
subgroup. Thus, the verbal closedness
(as well as algebraic closedness)
is a property of a
subgroup, while the strong verbal closedness is a property of an abstract
group.
The class of strongly verbally closed groups is fairly wide.
For example, the following
groups are strongly verbally closed:
\-
all abelian groups [Mazh18],
\-
all free groups [KM18],
\itemitem{--}
and even all virtually free groups
containing no
nonidentity finite normal subgroups
[KM18], [KMM18],
\-
all free products of
nontrivial groups [Mazh19],
\-
the fundamental group of any connected surface,
except for the Klein bottle [Mazh18], [K21],
\itemitem{--}
and even all acylindrically hyperbolic groups
without nonidentity finite normal subgroups [Bog23] (this generalises
several results above),
\-
all finite groups with nonabelian monoliths (in particular, all
finite simple groups) [KMO23],
\-
all finite symmetric groups [KK24],
\-
all dihedral groups of order not divisible by eight
[KMO23].

\enditem
There is also a general embedding theorem [KMO23]:
\disp{\sl
any group $H$ embeds into a strongly verbally closed group
of cardinality~$|H|+\aleph_0$
that satisfies all identities of $H$.
}%
Actually, neither proving nor  disproving the strong verbal closedness of
a given group is an easy task.

Until quite lately,
only few examples of
non-strongly-verbally-closed groups
were known.
However, the following theorem has recently been proven.

\Th KMO {\rm [KMO23]}.
The centre of a finite strongly verbally closed group
is its direct factor.

This theorem implies immediately
that, e.g., all finite nilpotent
non-abelian groups are  not
strongly verbally closed.
The same is true for finitely generated
nilpotent groups with non-abelian torsion subgroups [D23].
However, the following question from [KMO23] remains open:
\disp{\sl
does there exist a finitely generated
non-abelian nilpotent
strongly verbally closed group?
}%
(The embedding theorem [KMO23], see above, shows that the
finitely-generated condition cannot be dropped.)

The simplest infinite finitely generated
non-abelian nilpotent group (with abelian torsion subgroup) is the discrete
Heisenberg group $\UT_3(\Z)$, which turned out to be
non-strongly-verbally-closed [KMO23], but the proof of this fact is not
easy.

Let us ask a more daring question.

\Question.
Is it true that the centre of any finitely generated strongly verbally
closed group is its direct factor?

For now, we are able to prove only a weaker property --- the purity of
the centre,
and under some additional conditions (not
very onerous though). Recall that a subgroup $B$ of a group $A$ is called
\emph{pure} (or \emph{servant}) if,
for every integer $n$ and for every $b\in B$, the element $b$
is an $n$th power in $A$ if and only if it is an $n$th power in
$B$.

\proclaim Main theorem {\rm(a simplified form)}.
The centre of a finitely generated strongly verbally closed group
is pure if at least one of the following conditions holds\:
\-
all finite-order elements lie in the centre
{\rm(in this case, the purity of the centre is equivalent to its
isolatedness)},
\-
or the group is not infinitely decomposable.

\noindent
(The full statement can be found in the next section.)
Here, we call a group \emph{infinitely decomposable}
if, for any $n\in\N$, it decomposes into a product of $n$
pairwise commuting noncentral subgroups
(which are normal in this case of course).

Finitely generated infinitely decomposable groups are
exotic creatures, but they do exist; there exists even
a nontrivial finitely generated
group isomorphic to its Cartesian square [TJ74].

The central quotient of an infinitely decomposable group
decomposes into a direct product of arbitrarily many
nontrivial subgroups; this means that it does not satisfy the
maximality condition for direct factors by virtue of the following
theorem ([O95], Corollary 3 + Proposition 1):
\disp{\sl
a group with finitely generated
commutator quotient group
decomposes into a finite direct product of
indecomposable groups
if and only if it
satisfies the maximality condition for direct
factors;
and the number of such decompositions
\(up to
isomorphisms of the factors\)
is always finite.
}%
Thus, we obtain the following fact.

\Corollary 1.
The centre of a finitely generated strongly verbally closed group
is pure
if its central quotient
satisfies the maximality condition for direct factors
\(or, equivalently, decomposes into a finite direct product of
indecomposable groups\).
In particular, the centre of
any finitely generated solvable strongly verbally
closed group is pure.

This corollary of course contains the
mentioned above
fact that
the Klein-bottle group is not strongly verbally closed [K21]
(because the Klein-bottle group $K=\pres<a,b|a^2=b^2>$
is solvable, and its centre $\gp{a^2}$ is not pure:
$a^2$ is a square in $K$ but not in $\gp{a^2}$).
Moreover, our theorem gives many new interesting examples of
non-strongly-verbally-closed groups, e.g.,
\disp{\sl\hfuzz55pt
the following groups are not strongly verbally closed:
\itemitem{--}
all non-abelian braid groups
{\rm(because these groups are torsion-free,
and the centre is generated by the square of the
``fundamental" braid [G69])},
and, in particular, it is not 2-pure, i.e.,
not pure with respect to squares,
\itemitem{--}
$\SL_n(\Z)$ for even $n$;
{\rm indeed, for even $n$
the group $\SL_n(\Z)$ is not pure (precisely not 2-pure), because
its centre is generated by the square of the
block-diagonal matrix with blocks $\pmatrix{\phantom{-}0&1\cr-1&0}$,
and the group $\PSL_n(\Z)$ is indecomposable into a direct product,
as it is easy to show}.$^{*)}$
}%
\footnote{}{$^{*)}$
The braid groups are linear [Big01], [Kra02]
(as well as $\SL_n(\Z)$);
for linear (and even for all equationally
Noetherian) subgroups
$H$ of finitely generated groups $G$, the algebraic closedness
is equivalent to the existence of a retraction
(i.e., a homomorphism $G\to H$ that is
identical on $H$) [MR14];
hence, these examples can be stated in a ``more categorical"
language (see the abstract).
}%
The following corollary
(which we prove in the last section) is a
direct generalisation of Theorem KMO.

\Corollary 2.
If the commutator subgroup of a
finitely generated strongly verbally closed group is
finite \(or, equivalently, the index of the centre is finite%
\fnn{%
{\sl any group whose centre is of finite index
has finite commutator subgroup}
by Schur's theorem [Sch04].
The converse for finitely generated groups
was proven by B. Neumann [Neu51]
(but the history of this question is somewhat odd~[Ya16]).
}%
\), then the centre is a direct factor of the group.

Our approach to the proof of the main theorem may be called
similar to the proof from
[K21], but our argument is significantly trickier: in particular, we
use Chevalley's theorem on zeros of polynomials [Che35] (to be more
precise, we use its generalisation, see Section 3); in this aspect,
our argument is similar to the proof of Theorem KMO [KMO23].

{\noindent\bf The notation}
being used in the paper
is mainly standard. Note only that, if
$k\in \Z$, and $x$ and $y$ are elements of a group, then $x^y$,
$x^{ky}$, and $x^{-y}$ denote $y^{-1}xy$, $y^{-1}x^ky$,
and $y^{-1}x^{-1}y$,
respectively.
\emph{The commutator}~$[x,y]$ is $x^{-1}y^{-1}xy$.
The centre of a group $G$ is denoted by $\ZZ(G)$.
The cardinality of a set~$X$ is denoted by $|X|$.
If $X$ is a
subset of a group, then $\gp X$
denotes the subgroup generated by
$X$.
The symbol $\gp g_n$ denotes
the cyclic subgroup of order $n$
generated by an element $g$.
The letters~$\Z$~and~$\N$ stand for the sets of integers
and positive integers, respectively.
$\Z_n$ is the residues ring
modulo $n$ (or the additive group of this ring).
The free product of groups~$A$~and~$B$ is denoted by $A*B$.

The authors thank an anonymous referee
and A. A. Skutin
for useful comments allowing us to improve the readability of the paper.
The second author thanks
the Theoretical Physics and Mathematics Advancement Foundation ``BASIS".

\s 2.
Main theorem

First note that, to prove the main result as
it is stated in the introduction, we may restrict ourselves
to prime powers.

\disp{\sl
An abelian
subgroup $B$ of a group $A$ is pure \(= servant\) if and only if it
is
\emph{\hbox{$p$-pure} \(= $p$-servant\)} for each prime $p$, i.e.,
$B\cap\left\{a^{p^k}\;\bigm|\;a\in A\right\}=
\left\{b^{p^k}\;\bigm|\;b\in B\right\}$
for all positive integers $k$.
}%
Indeed, if $a^{mn}=b=b_1^m=b_2^n$, where
$m,n\in\Z$ are coprime,
$a\in A$, and $b,b_1,b_2\in B$,
then $(b_1^kb_2^l)^{mn}=b^{kn+lm}=b$,
where integers $k$ and $l$ are chosen such that $kn+lm=1$.


\proclaim Main theorem.
If a strongly verbally closed group $H$ is finitely generated
modulo its centre and
a central element~$b^{p^k}$ is not the $p^k$th power of any
central element \(where $b\in H$, the integer $p$ is prime, and $k\in\N$\),
then,
for any positive integer $n$,
the group $H$ decomposes into a product of $n$
pairwise commuting subgroups not
contained in
the
centraliser of the
set~$\left\{h\in H\;|\; h^{p^k}=1,\; [h,b]=1\right\}$.
{\rm(In particular, this set is not contained in the centre of $H$.)}

\Proof
This theorem follows immediately from the following
three lemmata, which we prove in the next section.

\proclaim Verbal-closedness lemma.
Suppose that $H$ is a group,
$b\in H$, and $b^{p^k}\in\ZZ(H)$, where $p$ is prime and $k\in\N$.
Then, in the semidirect product
$
G=F
\semitimes
\bigtimes\limits_{s\in S}
H_s,
$
where
\-
$S$ is a finite set,
\-
$F$ is a
subgroup in the group of functions $S\to\Z_{p^k}$
\(but $F$ is written multiplicatively\),
\-
$H_s$ are isomorphic copies of $H$
\(and the corresponding isomorphisms
are denoted by the subscripts:
$h\mapsto h_s$\),
\-
and
the action is the following:
$
h_s^f=\(h^{b^{f(s)}}\)_s,
$
\enditem
the
``diagonal" subgroup
$\left\{\prod\limits_{s\in S} h_s\;\Bigm|\; h\in H\right\}\iso H$
is verbally closed,
if each function $f\in F$ vanishes at some point~$s\in S$.

\proclaim Algebraic-non-closedness lemma.
Suppose that a group $H$ is finitely generated modulo its centre,
$b\in H$, and an element $b^{p^k}$, where $p$ is prime and $k\in\N$,
lies in the centre of $H$
but is not the
$p^k$th power of any central element.
Let $n$ be the maximal number
of pairwise commuting subgroups whose product is $H$
and none of which
centralises the set
$
\{h\in H\;|\; h^{p^k}=1,\; [h,b]=1\}
$.
Then
the diagonal subgroup $H$ of the semidirect product~$G$
from the
verbal-closedness lemma is not algebraically closed in $G$
if, for
any subset $T\subseteq S$ of cardinality at most
$\max(n,1)$,
there exists a function~$f_T\in F$ such that
$f_T(s)=1$ for all $s\in T$.

\noindent
To complete the proof of the theorem
(modulo these lemmata),
it remains to show
that the conditions of two lemmata above can be satisfied simultaneously
(for suitable $S$ and $F$); the following assertion
shows that they can.

\proclaim Function lemma.
For any positive integers $k$ and $n$ and any prime $p$, there exists
a finite set $S$ and a subgroup $F$ in the group
of functions $S\to\Z_{p^k}$ such that
\item{\rm1)}
each function from $F$ vanishes at some point of $S$;
\item{\rm2)}
but, for any $s_1,\dots,s_n\in S$, there exists a function $f\in F$ such
that $f(s_1)=\dots=f(s_n)=1$.

In the next section, we
prove these three lemmata, and thereby complete the
proof of the theorem.

\s 3.
Proofs of the three lemmata

{\noindent\bf Proof of the verbal-closedness lemma.}
To show that $H$ is verbally closed in $G$,
we have to prove that
any equation of the form
$$
w(x,y,\dots)=h,
\qbox{where $h\in H$ and
$w(x,y,\dots)$ is an element of a free group
$F(x,y,\dots)$,}
$$
which is solvable in $G$, is solvable in $H$.
A change of variables reduces
such an  equation
to the form
$$
x^mu(x,y,\dots)=h,
\qbox{where $h\in H$ and
$u(x,y,\dots)$ lies in the
commutator subgroup of a
free group
$F(x,y,\dots)$}.
$$
(To obtain this change of variables, note that
the free group modulo its commutator subgroup
is the free abelian group~$\Z\oplus\dots\oplus\Z$;
and, for any element $\^w$ of a free abelian group,
there is a basis $B$ such that $\^w=mb$,
where $m\in\Z$ and $b\in B$.)

\noindent
Suppose that this equation has a solution
$(\~x,\~y,\dots)$
in $G$;
in particular,
$
\~x=f^\epsilon d,
\hbox{where $\epsilon\in\Z$, $f\in F$,
and $d\in\!\bigtimes_{s\in S}\! H_s$}.
$
By the condition, the function $f$ vanishes at a point $s\in S$.
Consider the homomorphism (taking the $s$th coordinate):
$$
\phi\:G\to\~H=\gp\beta_{p^k}\semitimes H\
(\hbox{with the action }h^\beta=h^b),
\hbox{where
$\phi(h_s)=h$, $\phi(h_{s'})=1$ for $s'\ne s$,
$\phi(f')=\beta^{f'(s)}$
for $f'\in F$.
}
$$
Applying this homomorphism, we obtain that
the equation
$
w\bigl(\phi(\~x),y,z,\dots\bigr)=h
$
over $H$ with unknowns~$y,z,\dots$ (whose exponent sums
are zero in the word $w$) has a solution in the
group $\~H$.
(This is an equation over $H$, because $f(s)=0$ and, hence,
$\phi(\~x)\in H$.)
It remains to apply the following almost obvious fact:
\disp{\sl
if the exponent sum
of each variable in
a word $v\in H*F(y,z,\dots)$
is zero, and the equation $v=1$ is
solvable in the group $\~H$, then it is solvable in $H$ too.
}%
Indeed,
replacing $\beta$ with $b$ in a solution, we obtain a
solution from $H$.
This completes the proof.

\medskip

{\noindent\bf Proof of the algebraic-non-closedness lemma.}
Suppose that
$
H=\gp{h_1,\dots, h_m}\cdot\ZZ(H)
$.
The system of equations
$$
\left\{
\matrix{
x_f^{p^k}\=^1 1,\;
y_{i\,s}^{x_f}\=^2 y_{i\,s}^{b^{f(s)}},\;
\prod\limits_{q\in S}y_{i\,q}\=^3 h_i,\;
[y_{i\,s},y_{j\,s'}]\=^4 1,\;
[x_f,b]\=^5 1
}
\;\Biggm|\;f\in F,\;s\in S\ni s'\ne s,\;i,j=1,\dots,m\right\}
$$
over $H$ has an obvious solution in $G$
($
x_f=f,\;
y_{i\,s}=(h_i)_s
$)
but no solutions in $H$. Indeed,
\-
equations (4) mean that
the subgroups $K_s=\gp{y_{1\,s},\dots, y_{m\,s}}\cdot\ZZ(H)$
pairwise commute;
\-
equations (3) say that $H=\prod K_s$;
\-
the elements $x_f$ commute with $b$ by virtue of equation (5);
\-
all the aforesaid
observations and equation (1)
imply (by the condition of the lemma) that
the elements $x_f$ commute
with all
subgroups $K_s$, except for some $n$ (or less) subgroups~$K_t$,
where $t\in T\subseteq S$, $|T|\le n$
(this set $T$ is independent of $f$);

\-
thus, (by the condition of the lemma)
there exists a function $f_T\in F$ such that
$f_T(t)=1$ for all $t\in T$;
\-
hence,
equations (2) show that
the element $x_{f_T}b^{-1}$ commutes with $y_{i\,s}$ for $s\in T$;
\itemitem{--}
but, for $s\notin T$, the elements
$x_{f_T}b^{-1}$ and $y_{i\,s}$
commute too,
because, for such $s$, the element
$y_{i\,s}$ commutes with all $x_f$ (as mentioned above),
therefore, choosing a function $f\in F$ such that $f(s)=1$
(which exists by the condition), we obtain
$
y_{i\,s}=y_{i\,s}^{x_f}
\=^2
y_{i\,s}^{b^{f(s)}}=y_{i\,s}^b
$;
\-
thus, equations (3) show that the element
$x_{f_T}b^{-1}$ is central
(because $h_i$ generate $H$ modulo $\ZZ(H)$);
\-
raising this central element to the power $p^k$,
we obtain $b^{-p^k}$
(because
$(x_{f_T}b^{-1})^{p^k}
\=^5
x_{f_T}^{p^k}b^{-p^k}\=^1 b^{-p^k}$)
that contradicts the condition $b^{p^k}\notin\ZZ(H)^{p^k}$.
This completes the proof of the lemma.

\medskip

{\noindent\bf Proof of the function lemma.}
Put $S=\Z_{p^k}^m\setminus p\Z_{p^k}^m$, where $m$ is
sufficiently large (with respect to $k$
and $p$, see below).
Take as $F$ the set of (polynomial)
mappings~$\Z_{p^k}^m\to\Z_{p^k}$ defined by polynomials
(in $m$ variables)
without free terms and of degree at most $n(p-1)p^{k-1}$.

\noindent
Condition 2) holds. Indeed,
\-
we can choose a basis in
the free $\Z_{p^k}$-module
$\Z_{p^k}^m$ and assume
that all coordinates of
vectors $s_i$ vanish, except for the first $n$
coordinates (because the definition of $F$ does not depend on the choice
of a basis);
\-
consider the reduction modulo $p$
$\phi\:\Z_{p^k}^m\to\Z_{p^k}^m/p\Z_{p^k}^m\iso\Z_p^m$;
\-
each function from an $n$-dimensional
vector space
(containing
all $\phi(s_i)$) to  $\Z_p$ is represented by a polynomial whose
degree with respect to each variable is at most $p-1$
(because the polynomials $x^p$ and $x$ define the same function
$\Z_p\to\Z_p$);

\-
i.e., we have a polynomial $g\in\Z[x_1,\dots,x_n]$
of degree at most $n(p-1)$ such that
$g(0)\in p\Z_{p^k}$ and
$g(s_i)\in 1+p\Z_{p^k}$ for all $i$;
\-
hence, the polynomial $g^{p^{k-1}}$ takes (in $\Z_{p^k}$)
\itemitem{--}
value 0 at zero, because $p^{k-1}>k$
(so this
polynomial has no free term);
\itemitem{--}
value 1 at all points $s_i$ (because
$g(s_i)$ lie in $1+p\Z_{p^k}$, and this set is a (multiplicative) group
of
order $p^{k-1}$, so we raise an element of a group to the power
equal to the order of the group);
\-
thus,
$\deg g^{p^{k-1}}=p^{k-1}\cdot\deg g\le n(p-1)p^{k-1}$,
and $g^{p^{k-1}}\in F$ by the definition of $F$.

\enditem
To verify Condition 1), we apply a theorem
of Schanuel [Sch74], which says (in particular) that:
\disp{\sl
a polynomial
$f\in\Z_{p^k}[x_1,\dots,x_m]$ without free term
has
a nonzero root in $\{0,1\}^m$
if $m>(p^k-1)\cdot\deg f$.
}%
Thus, if we choose $m$ larger than $n(p-1)p^{k-1}(p^k-1)$,
then the both conditions are fulfilled. This completes the proof
of the lemma (and of the main theorem).

\s 4.
Proof of Corollary 2

This result follows easily from
the main theorem and
some known facts.
Indeed,
\-
the group $H$ is finitely generated,
and the index of the centre is finite, therefore, the centre
is a finitely generated abelian group
and, hence,
decomposes into a direct
product of a finite abelian group and a free abelian group $F$
(and the index $|H:F|$ is finite);
\-
as was shown by Schur [Sch04]
(see also [Rob96], Theorem 10.1.3),
\disp{\sl\hfuzz14.1pt
for any central finite-index subgroup $F$ of
\(any\) group~$H$,
the mapping
\newline
$\tau\:x\mapsto x^{|H:F|}$
\(called the \emph{transfer}\)
is a homomorphism from $H$ to $F$.
}%
\-
$\tau(H)=\tau(F)$,
because the centre is pure in $H$ (by the main theorem),
and $F$ is pure in the centre (obviously);
\-
therefore, the composition of $\tau$ and a (suitable)
isomorphism $\tau(H)=\tau(F)\to F$ is a retraction~$H\to F$;
\-
i.e., $F$ (being a central subgroup)
is a direct factor: $H=F\times H_1$;
\-
the second factor $H_1$ must be finite
(because $F$ is of finite index in $H$)
and strongly verbally closed:
\itemitem{--}
{\sl any direct factor of a strongly verbally closed group is
strongly verbally closed}%
\fn{%
The ``opposite" question of Mazhuga is still open:
{\sl is the class of strongly verbally closed groups
closed with respect to direct products?}
See [Mazh18d] for some results on this topic.}
(since both verbal and algebraic closedness
of a subgroup $A$ in a group
$\~A$
are equivalent to the corresponding properties
of~$A\times B$ in~$\~A\times B$);

\-
it remains to note that $\ZZ(H)=F\times\ZZ(H_1)$,
and the centre of $H_1$ is its direct factor by Theorem~KMO.

\References

[Big01]
S. Bigelow,
Braid groups are linear,
Journal of the American Mathematical Society, 14:2 (2001), 471-486.

[Bog22]
O. Bogopolski,
Equations in acylindrically hyperbolic groups and verbal closedness,
Groups, Geometry, and Dynamics, 16:2 (2022), 613-682.
\arXiv 1805.08071

[Che35]
C. Chevalley,
D\'emonstration d'une hypoth\`ese de M. Artin,
Abh. Math. Semin. Univ. Hambg., 11 (1935), 73-75.

[D23]
F. D. Denissov,
Finite normal subgroups of strongly verbally closed groups,
Journal of Group Theory (to appear).
\doi 10.1515/jgth-2023-0015
\arXiv:2301.02752

[G69]
F. A. Garside,
The braid group and other groups,
The Quarterly Journal of Mathematics, 20:1 (1969), 235-254.


[KK24]
O. K. Karimova, A. A. Klyachko,
Finite symmetric groups are strongly verbally closed,
arXiv:2405.01179\thinspace.

[K21]
A. A. Klyachko,
The Klein bottle group is not strongly verbally closed,
though awfully close to being so,
Canadian Mathematical Bulletin, 64:2 (2021), 491-497.
\arXiv 2006.15523

[KM18]
A. A. Klyachko, A. M. Mazhuga,
Verbally closed virtually free subgroups,
Sbornik: Mathematics, 209:6 (2018), 850-856.
\arXiv 1702.07761

[KMM18]
A. A. Klyachko, A. M. Mazhuga, V. Yu. Miroshnichenko,
Virtually free finite-normal-subgroup-free groups
are strongly verbally closed,
J. Algebra, 510 (2018), 319-330.
\arXiv 1712.03406

[KMO23]
A. A. Klyachko, V. Yu. Miroshnichenko, A. Yu. Olshanskii,
Finite and nilpotent strongly verbally closed groups,
Journal of Algebra and Its Applications 22:09 (2023), 2350188.
\arXiv:2109.12397

[Kra02]
D. Krammer,
Braid groups are linear,
Annals of Mathematics, 155:1 (2002), 131-156.
\arXiv math/0405198

[Mazh17]
A. M. Mazhuga,
On free decompositions of verbally closed subgroups
of free products of finite groups,
J.~Group Theory, 20:5 (2017), 971-986.
\arXiv 1605.01766

[Mazh18]
A. M. Mazhuga,
Strongly verbally closed groups,
J. Algebra, 493 (2018), 171-184.
\arXiv 1707.02464


[Mazh18d]
A. M. Mazhuga,
Verbally closed subgroups,
{\tencyrit Kandidat{s}kaya dissertaciya} ($\approx$ Ph.D. thesis),
MSU,
Moscow,
2018 (in Russian),
{\tt https://istina.msu.ru/dissertations/150563554/}\thinspace.

[Mazh19]
A. M. Mazhuga,
Free products of groups are strongly verbally closed,
Sbornik: Mathematics, 210:10 (2019), 1456-1492.
\arXiv 1803.10634

[MR14]
A. Myasnikov, V. Roman'kov,
Verbally closed subgroups of free groups,
J. Group Theory, 17:1 (2014), 29-40.
\arXiv 1201.0497

[Neu51]
B. H. Neumann,
Groups with finite classes of conjugate elements,
Proc. London Math. Soc., s3-1:1 (1951), 178-187.

[O95]
F. Oger,
The direct decompositions of a group $G$ with $G/G'$ finitely generated,
Trans. Amer. Math. Soc., 347:6 (1995), 1997-2010.

[Rob96]
D. J. S. Robinson,
A Course in the theory of groups.
Graduate texts in mathematics 80.
Springer, 1996.

[Rom12]
V. A. Roman'kov,
Equations over groups,
Groups - Complexity - Cryptology,
4:2 (2012), 191-239.

[RKh13]
V. A. Roman'kov, N. G. Khisamiev.
Verbally and existentially closed subgroups of free nilpotent groups.
Algebra and Logic, 52:4 (2013), 336-351.

[RKhK17]
V. A. Roman'kov, N. G. Khisamiev, A. A. Konyrkhanova,
Algebraically and verbally closed subgroups and retracts
of finitely generated nilpotent groups,
Siberian Math. J., 58:3 (2017), 536-545.

[RT20]
V. A. Roman'kov, E. I. Timoshenko,
Verbally Closed Subgroups of Free Solvable Groups,
Algebra and Logic, 59:3 (2020), 253-265.
\arXiv 1906.11689

[Sch74]
S. H. Schanuel,
An extension of Chevalley's theorem to congruences modulo prime powers,
Journal of Number Theory, 6:4 (1974), 284-290.

[Sch04]
J. Schur,
\"Uber die Darstellung der endlichen Gruppen durch gebrochene
lineare Substitutionen,
Journal f\"ur die reine und angewandte Mathematik, 127 (1904), 20-50.

[Tim21]
E. I. Timoshenko,
Retracts and verbally closed subgroups
with respect to relatively free soluble groups,
Sib. Math. J., 62:3, 537-544 (2021).

[TJ74]
J. M. Tyrer-Jones,
Direct products and the Hopf property,
J. Austral. Math. Soc., 17:2 (1974), 174-196.

[Ya16]
M. K. Yadav,
Central Quotient Versus Commutator Subgroup of Groups,
In: S. Rizvi, A. Ali, V. Filippis, (eds)
Algebra and its Applications.
Springer Proceedings in Mathematics \& Statistics, 174 (2016).
Springer, Singapore.
\arXiv 1011.2083

\end